\documentclass[10pt]{article}

\usepackage{amsmath}
\usepackage{amsfonts}
\usepackage{amsthm}
\usepackage[english]{babel}
\usepackage{graphicx}
\usepackage[all]{xy}
\usepackage{amssymb}
\usepackage{lineno}

\newtheorem{theorem}{Theorem}

\newtheorem{remark}{Remark}

\newcommand{\mR}{\mathbb{R}}

\begin{document}
\title{A new construction of the Clifford-Fourier kernel}

\author{Denis Constales$^{1}$ \footnote{E-mail: {\tt  denis.constales@ugent.be}} \and Hendrik De Bie$^{1}$\footnote{E-mail: {\tt hendrik.debie@ugent.be}} \and Pan Lian$^{1,2}$\footnote{E-mail: {\tt pan.lian@ugent.be}} }

\vspace{10mm}
\date{\small{1:  Department of Mathematical Analysis\\ Faculty of Engineering and Architecture -- Ghent University}\\
\small{Galglaan 2, 9000 Gent,
Belgium}\\ \vspace{5mm}
\small{2: Department of Mathematics - Harbin Institute of Technology\\
West Da-Zhi Street 92, 150001 Harbin, P.R.China}\\
\vspace{5mm}
}

\maketitle

\begin{abstract}
In this paper, we develop a new method based on the Laplace transform to study the Clifford-Fourier transform. First, the kernel of the Clifford-Fourier transform in the Laplace domain is obtained.
When the dimension is even, the inverse Laplace transform may be computed and we obtain the explicit expression for the kernel as a finite sum of Bessel functions.
We equally obtain the plane wave decomposition and find new integral representations for the kernel in all dimensions. Finally we define and compute the formal generating function for the even dimensional kernels.

~\\


\end{abstract}

\tableofcontents

\section{Introduction}
\label{sec:1}

In recent years, quite a bit of attention has been paid to the Clifford-Fourier transform
\[
F_{-}(f)(y)\int_{\mathbb{R}^{m}} K_m(x,y)f(x)dx
\]
with kernel given by
\[
K_{m}(x,y)=e^{ i\frac{\pi}{2}\Gamma_{y}}e^{-i(x,y)}
\]
where
\[
\Gamma_{y}:=-\sum_{j<k}e_{j}e_{k}(y_{j}\partial_{y_{k}}-y_{k}\partial_{y_{j}}).
\]
Here $(x,y)$ is the standard Euclidean inner product on $\mR^m$ and the $e_j$ generate a Clifford algebra.

This transform was introduced in \cite{FNF} as a generalization of the classical Fourier transform (FT) for multichannel signals. Because of the underlying Clifford algebra structure and its geometric interpretation, a transform is obtained that properly mixes the signals to be investigated. It nevertheless still satisfies many important properties of the classical FT, such as inversion, Plancherel theorem, behaviour of differentiation and (generalized) translation, etc.

It turned out to be a hard problem to compute its integral kernel $K_m$ explicitly. This was first done for $m=2$ in \cite{FNF2}. In \cite{FNF3} a recursive algorithm was presented to compute the kernel in even dimension. A completely explicit expression for all even dimensions was subsequently obtained in \cite{HY} by using plane wave decompositions.

In later work, the results were extended to fractional versions of the Clifford-Fourier transform \cite{MH, HN} and integral kernels satisfying certain generalized Helmholtz PDEs in Clifford analysis \cite{DBNS}. In \cite{Ro1} an approach using Lie superalgebras and group symmetries led to a complete classification of transforms that behave in the same way as both the Clifford-Fourier transform and the classical Fourier transform. Uncertainty principles for these transforms were obtained in \cite{GJ, Ro1}.

The main aim of the present paper is to develop a new and elegant method to compute the integral kernel $K_m$. This will be done by introducing an auxiliary variable $t$ and subsequently expressing its Laplace transform
\[
\mathcal{L}(t^{m/2-1} e^{-it(x,y)})
\]
in terms of the Cauchy kernel for the Dirac operator. In the Laplace domain, the action of $\Gamma_y$ is obtained using a monogenic expansion. Laplace inversion then yields our main result in the even dimensional case.

Our approach bypasses the need for a rather complicated induction argument in \cite{HY}. As an additional bonus, we are now able to compute an explicit generating function for all even dimensional kernels by again using the Laplace domain expression. This is achieved in Theorem \ref{th7} and \ref{th8}.

The paper is organized as follows. In Section 2 we recall basic facts concerning the Laplace transform, Clifford analysis and the Clifford-Fourier transform. In Section 3 we first compute the Laplace domain expression for the fractional Clifford-Fourier kernel. We use this result to reobtain both the plane wave decomposition and the explicit expression through Laplace inversion. Finally we derive a new integral identity for the kernel and we construct the generating function.

\section{Preliminaries}
\subsection{The Laplace transform}
The Laplace transform is widely used to solve differential and integral equations. Like the Fourier transform, the Laplace transform maps the original differential equation into an elementary algebraic expression. The solution of  the latter can then be transformed back to the solution of the original problem. Suppose that $f$ is a real or complex valued function of the variable $t>0$ and $s$ is a complex parameter. The Laplace transform of $f$ which has exponential order $\alpha$, i.e. $|f(t)|\le Ce^{\alpha t}, t\ge t_{0}$ is defined as
\[F(s)=\mathcal{L}(f(t))=\int_{0}^{\infty}e^{-st}f(t)dt.\] By Lerch's theorem \cite{Lp1}, if we restrict our attention to functions which are continuous on $[0,\infty)$, then the inverse transform \[\mathcal{L}^{-1}(F(s))=f(t)\] is uniquely defined. Inverse Laplace transforms can be computed directly by the complex inversion formula, which is based on  contour integration. Often, we can use integral transform tables (see e.g. \cite{lp2}) and the partial fraction expansion to compute Laplace transform. We list some  which will be used in this paper:
\begin{eqnarray}\mathcal{L}(t^{k-1})&=&\frac{\Gamma(k)}{s^{k}},  \quad k>0;\nonumber\\
\label{l7}\mathcal{L}(e^{-\alpha t})&=&\frac{1}{s+\alpha};\\
\label{l8} \mathcal{L}(t^{k-1}e^{-\alpha t})&=&\frac{\Gamma(k)}{(s+\alpha)^{k}},   \quad k>0;\\
 \label{l9} \mathcal{L}(\cos at)&= &\frac{s}{s^{2}+a^{2}};\\
\label{l10}\mathcal{L}(\sin at)&=&\frac{a}{ s^{2}+a^{2}}.
\end{eqnarray}
We also need some inverse Laplace transforms, with $r=(s^{2}+a^{2})^{1/2}$, $R=s+r$ and $g(s)=\mathcal{L}(f(t))$. We have
\begin{eqnarray}\label{l11} \mathcal{L}^{-1}(\frac{1}{r} (\frac{a}{R})^{\nu})&=&J_{\nu}(at),  \qquad \qquad \mbox{Re}(\nu) >-1, \mbox{Re} (s)>|\mbox{Im} (a)|;\\
\label{l2}\mathcal{L}^{-1}(2^{\nu}\pi^{-1/2}\Gamma(\nu+\frac{1}{2})a^{\nu}r^{-2\nu-1})&=&t^{\nu}J_{\nu}(at), \qquad   \mbox{Re}(\nu) >-1/2,  \mbox{Re} (s)>|\mbox{Im} (a)|;\\
\label{l3}\mathcal{L}^{-1}(2^{\nu+1}\pi^{-1/2}\Gamma(\nu+3/2)a^{\nu}r^{-2\nu-3}s )&=&t^{\nu+1}J_{\nu}(at),\qquad \mbox{Re}(\nu) >-1,  \mbox{Re} (s)>|\mbox{Im} (a)|; \\
\label{l4}\mathcal{L}^{-1}(r^{-1}g(r))&=&\int_{0}^{t}J_{0}[a(t^{2}-u^{2})^{1/2}]f(u)du;\\
\label{l5}\mathcal{L}^{-1}(g(r))&=&f(t)-a\int_{0}^{t}f[(t^{2}-u^{2})^{1/2}]J_{1}(au)du;\\
\label{l6}\mathcal{L}^{-1}(sr^{-1}g(r))&=&f(t)-at\int_{0}^{t}(t^{2}-u^{2})^{-1/2}J_{1}[a(t^{2}-u^{2})^{1/2}]f(u)du. \end{eqnarray}

The partial fraction decomposition of a rational polynomial
\begin{eqnarray*} F(s)=\frac{b_{m}s^{m}+b_{m-1}s^{m-1}+\cdots+b_{1}s+b_{0}}{a_{n}s^{n}+a_{n-1}s^{n-1}+\cdots+a_{1}s+a_{0}}=\frac{A(s)}{B(s)}, \qquad (n>m)
\end{eqnarray*}
expresses $F(s)$ as a sum of fractions with simple denominator. We only show the case when $F(s)$ has a single pole of order $m$. Then $F(s)$ can be expressed as
\begin{eqnarray*} F(s)=\frac{A(s)}{(s-p)^{m}}=\frac{c_{m}}{(s-p)^{m}}+\frac{c_{m-1}}{(s-p)^{m-1}}+\cdots+\frac{c_{1}}{s-p}
,\end{eqnarray*}
with complex constants $c_{m-j}=\frac{1}{j!}\frac{d^{m-j}}{d s^{m-j}}[F(s)(s-p)^{m}]_{s=p}, j=1,\cdots, m.$

The Laplace transform of a matrix valued function is simply the matrix of Laplace transforms of the individual elements. For example \[\mathcal{L}\begin{pmatrix} e^{t}\\te^{-t}\end{pmatrix} =\begin{pmatrix} 1/(s-1)\\1/(s+1)^{2}\end{pmatrix}.\]
Suppose $A$ is an $n\times n$ matrix. The matrix exponential is interpreted in terms of a power series, namely\[\exp(At)=I +At+\frac{A^{2}t^{2}}{2!}+\frac{A^{3}t^{3}}{3!}+\ldots.\]
By analogy with the scalar case, we have
\[\mathcal{L}(e^{At})=(sI-A)^{-1}.\] For more about the Laplace transform of matrix-valued functions, see \cite{M}.

\subsection{Clifford analysis }
Let $\mathbb{R}^{m}$ be the usual $m$-dimensional Euclidean space with an orthonormal basis $\{e_{1}, e_{2}, \ldots, e_{m}\}$. The Clifford algebra $\mathcal{C}\ell_{0,m}$ associated with $\mathbb{R}^{m}$ is spanned by the reduced products \[\mathop{\cup}_{j=1}^{m}\{e_{\alpha}=e_{i_{1}}e_{i_{2}}\ldots e_{i_{j}}:\alpha=\{i_{1},i_{2},\ldots, i_{j}\},  \quad  1\le i_{1}<i_{2}<\cdots<i_{j}\le m\}\]  with the relations $e_{i}e_{j}+e_{j}e_{i}=-2\delta_{ij}$. We have $\mathcal{C}\ell_{0,m}=\{\sum_{\alpha}e_{\alpha}x_{\alpha}; x_{\alpha}\in \mathbb{R} \}$.
The Clifford algebra $\mathcal{C}\ell_{0,m}$ is a graded algebra as  $\mathcal{C}\ell_{0,m}=\oplus_{l}\mathcal{C}\ell_{0,m}^{l}$ where $\mathcal{C}\ell_{0,m}^{l}$ is spanned by reduced Clifford products of length $l$. The Clifford algebra $\mathcal{C}\ell_{0,m}$ is a $2^{m}$-dimensional real associative algebra with identity and contains a copy of $\mathbb{R}^{m}$ by the canonical mapping $x=\sum_{j=1}^{n} x_{j}e_{i}$. Now we can define the inner product and the wedge product of two vectors $x,y\in \mathbb{R}^{m}$ using the Clifford product: \begin{eqnarray*}(x,y):&=&\sum_{j=1}^{m}x_{j}y_{j}=-\frac{1}{2}(xy+yx);\\
x\wedge y:&=&\sum_{j<k}e_{j}e_{k}(x_{j}y_{k}-x_{k}y_{j})=\frac{1}{2}(xy-yx).\end{eqnarray*}
It is easy to get $xy=-(x,y)+x\wedge y$, and $(x\wedge y)^{2}=-|x|^{2}|y|^{2}+(x,y)^{2}$ (see \cite{HY}).
Because $\frac{(x\wedge y)^{2}}{|x\wedge y|^{2}}=-1,$ we can consider $\frac{x\wedge y}{|x\wedge y|}$ as an imaginary unit. Now we can Laplace transform $e^{-(x\wedge y) t}$ and get \begin{eqnarray*}\mathcal{L}(e^{-(x\wedge y) t})=\frac{1}{s+x\wedge y}.\end{eqnarray*}
The complexified Clifford algebra $\mathcal{C}\ell_{0,m}^{c}$ is defined as $\mathbb{C}\otimes \mathcal{C}\ell_{0,m}$.

The conjugation is defined by $\overline{(e_{j_{1}}\ldots e_{j_{l}})}= (-1)^{l}e_{j_{l}}\ldots e_{j_{1}}$ as a linear mapping. For $x,y\in \mathcal{C}\ell_{0,m}^{c}$, we have $\overline{(xy)}=\overline{y}\overline{x}, \overline{\overline{x}}=x,$ and $\overline{i}=i$ which is not the usual complex conjugation. We define the Clifford norm of $x$ by $|x|^{2}=x\overline{x}, x\in \mathcal{C}\ell_{0,m}^{c}.$
\begin{remark} Even when we complexify  $x$, we shall use the sum of squares norm, not the hermitian norm.
\end{remark}

The Dirac operator is defined as:
\[D=\sum_{j=1}^{m}e_{j}\partial_{x_{j}}.\]  When $u$ is a scalar $C^{1}$ function, $Du$ can be identified with the gradient $\nabla u$. A function is called monogenic if $Du=0$. An important example of a monogenic function is the generalized Cauchy kernel \cite{ca}
\[G(x)=\frac{1}{\omega_{m}}\frac{\bar{x}}{|x|^{m}}\] where $\omega_{m}$ is the surface area of the unit ball in $\mathbb{R}^{m}$. It is the fundamental solution of the Dirac operator. Note that the norm here is $|x|=(\sum_{i=1}^{m}x_{i}^{2})^{1/2}$ and coincides with the Clifford norm. In the following, when the two norms are equal, we will not point it out again.

Denote by $\mathcal{P}$ the space of polynomials taking values in $\mathcal{C}\ell_{0,m}$, i.e. $\mathcal{P}:=\mathbb{R}[x_{1},\ldots,x_{m}]\otimes\mathcal{C}\ell_{0,m}$. The space of homogeneous polynomials of degree $k$ is then denoted by $\mathcal{P}_{k}$. The space $\mathcal{M}_{k}:=(\mbox{ker}D)\cap \mathcal{P}_{k}$  is called the space of spherical monogenics of degree $k$.

The local behaviour of a monogenic function near a point can be investigated by the polynomials introduced above. The following theorem  is  the analogue of the Taylor series in complex analysis.
\begin{theorem} \label{ts}\cite{ca} Suppose $f$ is monogenic in an open set $\Omega$ containing the origin.  Then there exists an open neighbourhood $\Lambda$ of the origin in which $f$ can be developed into a normally convergent series of spherical monogenics $ M_{k}f(x)$, i.e.
\[f(x)=\sum_{k=0}^{\infty}M_{k}f(x),\]
with $M_{k}f(x)\in \mathcal{M}_{k}$.
\end{theorem}

We further introduce the Gamma operator (see \cite{ca})
\[\Gamma_{x}:=-\sum_{j<k}e_{j}e_{k}(x_{j}\partial_{x_{k}}-x_{k}\partial_{x_{j}})=-x D_{x}-E_{x}.\]
Here $E_{x}=\sum_{i=1}^{m}x_{i}\partial_{x_{i}}$ is the Euler operator.  Note that $\Gamma_{x} $ commutes with scalar radial functions. The operator $\Gamma_{x}$ has two important eigenspaces:
\begin{eqnarray}\label{eg}\Gamma_{x}\mathcal{M}_{k}=-k\mathcal{M}_{k},\end{eqnarray}\begin{eqnarray}\label{eg1}\Gamma_{x}(x\mathcal{M}_{k-1})=(k+m-2)x\mathcal{M}_{k-1}\end{eqnarray}
which follows from the definition of $\Gamma_{x}$.
\subsection{The Clifford-Fourier transform}
The classical Fourier transform \[\mathcal{F}(f)(y)=(2\pi)^{-m/2}\int_{\mathbb{R}^{m}}e^{-i(x,y)}f(x)dx,\]
with $(x,y)$  the inner product on $\mathbb{R}^{m}$, can be represented by the operator exponential \[\mathcal{F}=e^{\frac{i\pi m}{4}}e^{\frac{i\pi}{4}(\Delta-|x|^{2})},\] see e.g. \cite{HR}, \cite{FG}.
Brackx, De Schepper and Sommen introduced a pair of  Fourier transforms using the angular Dirac operator $\Gamma_{x}$ in the Clifford algebra setting in \cite{FNF}. More precisely, it is defined by \[\mathcal{F}_{\pm}=e^{\frac{i\pi m}{4}}e^{\mp\frac{i\pi}{2}\Gamma_{x}}e^{\frac{i\pi}{4}(\Delta-|x|^{2})}.\]
For $F_{-}$, we denote the kernel as
\[K_{ m}(x,y)=e^{ i\frac{\pi}{2}\Gamma_{y}}e^{-i(x,y)}. \] In general, it is not easy to compute this kernel explicitly. In \cite{HY}, the authors derived the kernel for  even dimensions as a finite sum of Bessel functions. Later in \cite{HN}, the fractional Clifford-Fourier transform was introduced as a generalization of the fractional Fourier transform \[\mathcal{F}_{\alpha,\beta}=e^{\frac{i\alpha m}{2}}e^{i\beta\Gamma_{x}}e^{\frac{i\alpha}{2}(\Delta-|x|^{2})}\]
and the kernels of even dimensions were obtained by a similar method. In \cite{MH},  a new construction of the fractional Clifford-Fourier kernels was given by solving wave-type problems. In the present paper the fractional Clifford-Fourier kernel is computed as \[K_{m}^{p}(x,y)=e^{ip\Gamma_{y}}e^{-i(x,y)}.\] The more general case can also be obtained using our method.

\section{Laplace transform method}
\subsection{The fractional Clifford-Fourier kernel in the Laplace domain}
In this section we introduce an auxiliary variable $t$ in the exponent of the classical Fourier transform and then use the Laplace transform to get the Clifford-Fourier kernel in the Laplace domain.

 We use the notation $\sqrt{-}:=\sqrt{s^{2}-|x|^{2}|y|^{2}}$. By direct computation, $(s+\sqrt{-})(s-\sqrt{-})=|x|^{2}|y|^{2}$. We also have
\begin{eqnarray}\label{f1}
\left|1+\frac{yx}{s+\sqrt{-}}\right|^{2}&=&\biggl(1+\frac{yx}{s+\sqrt{-}}\biggr)\overline{\biggl(1+\frac{yx}{s+\sqrt{-}}\biggr)}\nonumber\\&=&\biggl(1+\frac{yx}{s+\sqrt{-}}\biggr)\biggl(1+\frac{xy}{s+\sqrt{-}}\biggr)\nonumber\\
&=&1+\frac{yx+xy}{s+\sqrt{-}}+\frac{|x|^{2}|y|^{2}}{(s+\sqrt{-})^{2}}\nonumber\\&=&1-\frac{2(x,y)}{s+\sqrt{-}}+\frac{(s+\sqrt{-})(s-\sqrt{-})}{(s+\sqrt{-})^{2}}\nonumber\\&=&\frac{2(s-(x,y))}{s+\sqrt{-}},
\end{eqnarray}
where we have used the Clifford norm.
Then using (\ref{f1}), we can express $\mathcal{L}(t^{m/2-1}e^{t(x,y)})$ in terms of the generalized Cauchy kernel introduced in the previous section. We have
\begin{eqnarray} \label{lp1} \mathcal{L}(\Gamma(m/2)e^{t(x,y)})&=&\frac{\Gamma(m/2)}{(s-(x,y))^{m/2}}\nonumber\\&=&\frac{\Gamma(m/2)}{(\displaystyle \frac{s+\sqrt{-}}{2})^{m/2}\left|1+\frac{yx}{s+\sqrt{-}}\right|^{m}}\nonumber\\&=&
\frac{\Gamma(m/2)}{(\frac{s+\sqrt{-}}{2})^{m/2}\left|1+\frac{yx}{s+\sqrt{-}}\right|^{m}}\frac{\displaystyle 1+\frac{yx}{s+\sqrt{-}}-\frac{\displaystyle y(1+\frac{yx}{s+\sqrt{-}})x}{s+\sqrt{-}}}{\displaystyle \frac{2\sqrt{-}}{\displaystyle s+\sqrt{-}}}\nonumber\\&=&
\frac{2^{m/2-1}\Gamma(m/2)}{\sqrt{-}(s+\sqrt{-})^{m/2-1}}
\frac{\displaystyle 1+\frac{yx}{s+\sqrt{-}}-\frac{ y(1+\frac{ \displaystyle yx}{\displaystyle s+\sqrt{-}})x}{s+\sqrt{-}}}{\displaystyle \left|1+\frac{yx}{s+\sqrt{-}}\right|^{m}}.
\end{eqnarray}

The first equality is by (\ref{l8}),  the second equality by (\ref{f1}), and the third equality follows by
\begin{eqnarray*} 1+\frac{yx}{s+\sqrt{-}}-\frac{y(1+\frac{yx}{s+\sqrt{-}})x}{s+\sqrt{-}}&=&1+\frac{yx}{s+\sqrt{-}}-\frac{yx}{s+\sqrt{-}}-\frac{yyxx}{(s+\sqrt{-})^{2}}\\&=&1-\frac{(s+\sqrt{-})(s-\sqrt{-})}{(s+\sqrt{-})^{2}}\\&=&\frac{2\sqrt{-}}{s+\sqrt{-}}.
\end{eqnarray*}

Next we will compute $\mathcal{L}(t^{m/2-1}e^{ip\Gamma_{y}}e^{t(x,y)})$ by acting with $e^{ip\Gamma_{y}}$ on both sides of (\ref{lp1}).
The generalized Cauchy kernel $G(y)=\frac{1}{\omega_{m}}\frac{\bar{y}}{|y|^{m}}$ is a monogenic function except at the origin. By translation,  $\frac{y+x}{|y+x|^{m}}$, $\frac{y+\frac{1}{x}}{|y+\frac{1}{x}|^{m}}\frac{x}{|x|^{m}}$, $\frac{yx+1}{|yx+1|^{m}}$ are also monogenic in $y$ except at $-x,-x^{-1},-x^{-1}$ respectively. Using Theorem \ref{ts}, we can express $\frac{yx+1}{|yx+1|^{m}}$ as a series of spherical monogenic polynomials, i.e. \[\frac{yx+1}{|yx+1|^{m}}=M_{0}(y)+M_{1}(y)+M_{2}(y)+\cdots\] where $M_{k}(y)$ is a  spherical monogenic of order $k$.
Substituting $\frac{y}{s+\sqrt{-}}$ for $y$, we have \[ \frac{\displaystyle 1+\frac{yx}{s+\sqrt{-}}}{\left|\displaystyle 1+\frac{yx}{s+\sqrt{-}}\right|^{m}}=\frac{M_{0}(y)}{(s+\sqrt{-})^{0}}+\frac{M_{1}(y)}{(s+\sqrt{-})^{1}}+\cdots.\]
Using (\ref{eg}), we obtain
\begin{eqnarray*} &&\Gamma_{y}\biggl(\frac{\displaystyle 1+\frac{yx}{s+\sqrt{-}}}{\left|\displaystyle 1+\frac{yx}{s+\sqrt{-}}\right|^{m}}\biggr)\\&=&\frac{\Gamma_{y} M_{0}(y)}{(s+\sqrt{-})^{0}}+\frac{\Gamma_{y} M_{1}(y)}{(s+\sqrt{-})^{1}}+\cdots\\&=&\frac{0\cdot M_{0}(y)}{(s+\sqrt{-})^{0}}+\frac{(-1)\cdot M_{1}(y)}{(s+\sqrt{-})^{1}}+\frac{(-2)\cdot M_{2}(y)}{(s+\sqrt{-})^{2}}+\cdots\end{eqnarray*}
and so
\begin{eqnarray*}&&e^{ip\Gamma_{y}}\biggl(\frac{\displaystyle 1+\frac{yx}{s+\sqrt{-}}}{\left|\displaystyle 1+\frac{yx}{s+\sqrt{-}}\right|^{m}}\biggr)\\&=&\frac{e^{ip\cdot 0}M_{0}(y)}{(s+\sqrt{-})^{0}}+\frac{e^{ip\cdot (-1)}M_{1}(y)}{(s+\sqrt{-})^{1}}+\frac{e^{ip\cdot (-2)}M_{2}(y)}{(s+\sqrt{-})^{2}}+\cdots\\&=&M_{0}\biggl(\frac{e^{-ip}y}{(s+\sqrt{-})}\biggr)+M_{1}\biggl(\frac{e^{-ip}y}{(s+\sqrt{-})}\biggr)+
M_{2}\biggl(\frac{e^{-ip}y}{(s+\sqrt{-})}\biggr)+\cdots\\&=& \frac{\displaystyle 1+\frac{e^{-ip}yx}{s+\sqrt{-}}}{\displaystyle \left|1+\frac{e^{-ip}yx}{s+\sqrt{-}}\right|^{m}}.
\end{eqnarray*}
Similarly, by (\ref{eg1}),  \[e^{ip\Gamma_{y}}(yM_{k}(y))=e^{i(m-2)p+i(k+1)p}(yM_{k}(y))=e^{i(m-2)p}(ye^{ip}M_{k}(e^{ip}y)).\] Now we can get the desired result \begin{eqnarray*}
&&\mathcal{L}(t^{m/2-1}e^{ip\Gamma_{y}}e^{t(x,y)})\\&=&\frac{2^{m/2-1}(m/2-1)!}{\sqrt{-}(s+\sqrt{-})^{m/2-1}}\biggl(\frac{1+\frac{\displaystyle e^{-ip}yx}{\displaystyle s+\sqrt{-}}}{\displaystyle \left|1+\frac{e^{-ip}yx}{s+\sqrt{-}}\right|^{m}}-e^{i(m-2)p}\frac{\frac{\displaystyle e^{ip}y(1+\frac{e^{ip}yx}{\displaystyle s+\sqrt{-}})x}{s+\sqrt{-}}}{\left|1+\frac{\displaystyle e^{ip}yx}{\displaystyle s+\sqrt{-}}\right|^{m}}\biggr).
\end{eqnarray*}
In order to simplify the expression further, we need the following, \begin{eqnarray} \label{f2} \left|1+\frac{e^{-ip}yx}{s+\sqrt{-}}\right|^{m}&=&[(1+\frac{e^{-ip}yx}{s+\sqrt{-}})(1+\frac{e^{-ip}xy}{s+\sqrt{-}})]^{m/2}\nonumber\\&=&[1-\frac{2e^{-ip}(x,y)}{s+\sqrt{-}}+\frac{e^{-2ip}yxxy(s-\sqrt{-})}{(s+\sqrt{-})^{2}(s-\sqrt{-})}]^{m/2}\nonumber\\&=&
[\frac{2e^{-ip}}{s+\sqrt{-}}(1/2e^{ip}(s+\sqrt{-})-(x,y)+1/2e^{-ip}(s-\sqrt{-}))]^{\frac{m}{2}}
\nonumber\\&=&[\frac{2e^{-ip}}{s+\sqrt{-}}(s\cos p-(x,y)+i\sqrt{-}\sin p)]^{\frac{m}{2}},
\end{eqnarray}
as well as
\begin{eqnarray}\label{f3} 1+\frac{e^{-ip}yx}{s+\sqrt{-}}=\frac{s+\sqrt{-}+e^{-ip}yx}{s+\sqrt{-}}
\end{eqnarray}
and
\begin{eqnarray}\label{f4} \frac{y(1+\frac{e^{ip}yx}{s+\sqrt{-}})x}{s+\sqrt{-}}=\frac{yx+e^{ip}(s-\sqrt{-})}{s+\sqrt{-}}=\frac{e^{ip}(e^{-ip}yx+s-\sqrt{-})}{s+\sqrt{-}}
.\end{eqnarray}
Combining (\ref{f2}), (\ref{f3}) and (\ref{f4}), we get the following theorem.
\begin{theorem}\label{qq} The Laplace transform of the fractional Clifford-Fourier kernel is given by:
 \begin{eqnarray*}\mathcal{L}(t^{m/2-1}e^{ip\Gamma_{y}}e^{t(x,y)})=\frac{\Gamma(m/2)}{2\sqrt{-}}\biggl(\frac{s+\sqrt{-}+e^{-ip}yx}{(e^{-ip}(s\cos p+i\sqrt{-}\sin p-(x,y)))^{m/2}}\nonumber\\-e^{imp}\frac{s-\sqrt{-}+e^{-ip}yx}{(e^{ip}(s\cos p-i\sqrt{-}\sin p-(x,y)))^{m/2}}\biggr),
\end{eqnarray*}
with $\sqrt{-}=\sqrt{s^{2}-|x|^{2}|y|^{2}}$.
\end{theorem}
If we substitute $-ix$ for $x$ in Theorem \ref{qq}, and denote $\sqrt{+}:=\sqrt{s^{2}+|x|^{2}|y|^{2}}$, note that the norm here is the Clifford norm, we have the following theorem which can be used to get the kernel by transforming back and setting $t=1$.
\begin{theorem}\label{qq1} The Laplace transform of the fractional Clifford-Fourier kernel is given by:
 \begin{eqnarray*}\mathcal{L}(t^{m/2-1}e^{ip\Gamma_{y}}e^{-it(x,y)})=\frac{\Gamma(m/2)}{2\sqrt{+}}\biggl(\frac{s+\sqrt{+}-ie^{-ip}yx}{(e^{-ip}(s\cos p+i\sqrt{+}\sin p+i(x,y)))^{m/2}}\nonumber\\-e^{imp}\frac{s-\sqrt{+}-ie^{-ip}yx}{(e^{ip}(s\cos p-i\sqrt{+}\sin p+i(x,y)))^{m/2}}\biggr),
\end{eqnarray*}
with $\sqrt{+}=\sqrt{s^{2}+|x|^{2}|y|^{2}}$.
\end{theorem}
\subsection{Plane wave decompostion of the fractional Clifford-Fourier kernel}
In this subsection we use the notation $\hat{x}=\frac{x}{|x|}$, $\hat{y}=\frac{y}{|y|}$ to denote two unit vectors. For $\hat{x}$, $\hat{y}$ we also have the result in Theorem \ref{qq1}. This time, we could get the kernel by putting $t=|x||y|$. Denote $r=\sqrt{s^{2}+1}$, $R=s+\sqrt{s^{2}+1}$, and $(\hat{x},\hat{y})=\cos \theta$. Using $s=\frac{R-1/R}{2}$ and $\sqrt{s^{2}+1}=\frac{R+1/R}{2}$, Theorem \ref{qq1} becomes
\begin{eqnarray}\label{add1}&&\mathcal{L}(t^{m/2-1}e^{ip\Gamma_{y}}e^{-it(\hat{x},\hat{y})})\nonumber\\&=&\frac{\Gamma(m/2)}{r}2^{m/2-1}R^{-m/2}\biggl(\frac{R-ie^{-ip}\hat{y}\hat{x}}{(1+2\frac{ie^{-ip}}{R}\cos \theta -(\frac{e^{-ip}}{R})^{2})^{m/2}}+e^{imp}\frac{\frac{1}{R}+ie^{-ip}\hat{y}\hat{x}}{(1+2\frac{ie^{ip}}{R}\cos \theta -(\frac{e^{ip}}{R})^{2})^{m/2}}\biggr)
\nonumber\\&=&\frac{\Gamma(m/2)}{r}2^{m/2-1}R^{-m/2}\nonumber\\&& \times\biggl(\frac{1}{(1+2\frac{ie^{-ip}}{R}\cos \theta -(\frac{e^{-ip}}{R})^{2})^{m/2}}[ie^{-ip}(-\cos \theta + \frac{-ie^{-ip}}{R})+(R+2ie^{-ip}\cos \theta-\frac{e^{-2ip}}{R})+ie^{-ip}\hat{x}\wedge\hat{y}]\nonumber\\&&+e^{imp}\frac{\frac{1}{R}+ie^{-ip}( -\cos\theta+\hat{y}\wedge\hat{x})}{(1+2\frac{ie^{ip}}{R}\cos \theta -(\frac{e^{ip}}{R})^{2})^{m/2}}\biggr) \end{eqnarray}
Recall the generating function of the Gegenbauer polynomial \cite{sw} \begin{eqnarray}\label{geg2}\frac{1}{(1-2xt+t^{2})^{\lambda}}=\sum_{k=0}^{\infty}C^{(\lambda)}_{k}(x)t^{k},\end{eqnarray}  and its  derivative with respect to $t$ \begin{eqnarray}\label{geg1}-\lambda\frac{-2x+2t}{(1-2xt+t^{2})^{\lambda+1}}=\sum_{k=0}^{\infty}kC^{(\lambda)}_{k}(x)t^{k-1}.\end{eqnarray}

Note that (\ref{add1}) consists of  five parts. We use (\ref{geg2})(e.g. the wedge term) and (\ref{geg1})(e.g. second term) to write each of them as a series. When transforming back by (\ref{l11}),
 we get the plane wave decomposition of the fractional Clifford-Fourier kernel as follows which can be compared with Theorem 3.2 in \cite{HN}.
\begin{theorem}  The series representation of the fractional Clifford-Fourier kernel is given by
\begin{eqnarray*}K_{m}^{p}(x,y)&=&e^{ip\Gamma_{y}}e^{-i(x,y)}\\&=&A^{p}_{m}+B^{p}_{m}+x\wedge yC^{p}_{m},
\end{eqnarray*}
where \begin{eqnarray*} A^{p}_{m}&=&-2^{m/2-2}\Gamma(m/2)\sum_{k=0}^{\infty}i^{-k}(e^{ip(k+m-2)}-e^{-ipk})(|x||y|)^{-m/2+1}J_{m/2-1+k}(|x||y|)C^{(m/2-1)}_{k}(\cos\theta),\\
B^{p}_{m}&=&2^{m/2-2}\Gamma(m/2)\sum_{k=0}^{\infty}i^{-k}(k+m/2-1)(e^{ip(k+m-2)}+e^{-ipk})(|x||y|)^{-m/2+1}J_{m/2-1+k}(|x||y|)C^{(m/2-1)}_{k}(\cos\theta),
\\C^{p}_{m}&=&2^{m/2-1}\Gamma(m/2) \sum_{k=1}^{\infty}i^{-k}(e^{ip(k+m-2)}-e^{-ipk})(|x||y|)^{-m/2}J_{m/2-1+k}(|x||y|)C^{(m/2)}_{k-1}(\cos\theta).
\end{eqnarray*}
\end{theorem}
Alternatively,  using the generating function of  Gegenbauer polynomials, we have
\begin{eqnarray*}&&(1+2\cos \theta\frac{ie^{-ip}}{R}+(\frac{-ie^{-ip}}{R})^{2})^{-m/2}\nonumber\\&=&\sum_{k=0}^{\infty}(\frac{-ie^{-ip}}{R})^{k} C_{k}^{(m/2)}(\cos \theta)\nonumber\\&=&\sum_{k=0}^{\infty} \sum_{a=0}^{k}(\frac{-ie^{-ip}}{R})^{k}\frac{(m/2)_{a}(m/2)_{k-a}}{a!(k-a)!}\cos(k-2a)\theta\end{eqnarray*}
which means we can express  formula (\ref{add1}) equally as a Fourier series.

\subsection{Even dimensional Clifford-Fourier kernel}

When $p=\pi/2$, the  result in Theorem \ref{qq1} reduces to
\begin{eqnarray}\label{f5}&&\mathcal{L}(t^{m/2-1}e^{i\frac{\pi}{2}\Gamma_{y}}e^{-it(x,y)})\nonumber\\&=&\frac{\Gamma(m/2)}{2\sqrt{+}}\biggl(\frac{s+\sqrt{+}-yx}{(\sqrt{+}+(x,y))^{m/2}}-e^{im\pi/2}\frac{s-\sqrt{+}-yx}{(\sqrt{+}-(x,y))^{m/2}}\biggr)
\nonumber\\&=&\frac{\Gamma(m/2)}{2\sqrt{+}}
\nonumber\\&&\biggl(\frac{(s-yx+\sqrt{+})(\sqrt{+}-(x,y))^{m/2}-e^{im\pi/2}(\sqrt{+}+(x,y))^{m/2}(s-yx-\sqrt{+})}{(s^{2}+(ix\wedge y)^{2})^{m/2}}\biggl).\nonumber\\\end{eqnarray}
When $m/2$ is even, (\ref{f5}) becomes \begin{eqnarray}\label{ab1}&&(m/2-1)!\biggl(\frac{(s-yx)(\sum_{j=1,3,5,\cdots}\binom{\frac{m}{2}}{j}(\sqrt{+})^{m/2-j-1}(-(x,y))^{j})}{(s^{2}+(ix\wedge y)^{2})^{m/2}}\nonumber\\&&+\frac{\sum_{j=0,2,4,\cdots}\binom{\frac{m}{2}}{j}(\sqrt{+})^{m/2-j}(-(x,y))^{j}}{(s^{2}+(ix\wedge y)^{2})^{m/2}}\biggr)\nonumber\\&=&(m/2-1)!\nonumber\\&&\biggl(\frac{(s-yx)(\sum_{j=1,3,5,\cdots}\binom{\frac{m}{2}}{j}(s^{2}+(ix\wedge y)^{2}-(x,y)^{2})^{\frac{m/2-j-1}{2}}(-(x,y))^{j})}{(s^{2}+|x\wedge y|^{2})^{m/2}}\nonumber\\&&+\frac{\sum_{j=0,2,4,\cdots}\binom{\frac{m}{2}}{j}(\sqrt{+})^{m/2-j}(-(x,y))^{j}}{(s^{2}+|x\wedge y|^{2})^{m/2}}\biggr)
\end{eqnarray} where all the sums are finite.
Simplifying (\ref{ab1}), we find that it is a finite sum of polynomials of type $\frac{(x,y)^{k}}{(s^{2}+|x\wedge y|^{2})^{q}}$ and  $\frac{s (x,y)^{k}}{(s^{2}+|x\wedge y|^{2})^{q}}$. Formulas (\ref{l2}) and (\ref{l3}) show that the kernel can be expressed as a finite sum of Bessel functions. Now we can get the kernel expressed in terms of Bessel functions which has been obtained in a completely different way in \cite{HY}.
\begin{theorem} The kernel of the Clifford-Fourier transform for even dimension $m=4n, n\ge 1$ is given by \begin{eqnarray*}&&K_{m}(x,y)=e^{i\frac{\pi}{2}\Gamma_{y}}e^{-i(x,y)}\\&=&(\pi/2)^{1/2}\biggl(A_{m}(u,v)+B_{m}(u,v)+(x\wedge y)C_{m}(u,v)\biggr) \end{eqnarray*} where $u=(x,y)$ and $v=|x\wedge y|$ and
 \begin{eqnarray*} A_{m}(u,v)=\sum^{m/4-1}_{l=0}u^{m/2-2-2l}\frac{1}{2^{l}l!}\frac{(m/2)!}{(m/2-2l-1)!}\frac{J_{(m-2l-3)/2}(v)}{v^{(m-2l-3)/2}},\end{eqnarray*}
 \begin{eqnarray*}
B_{m}(u,v)=-\sum^{m/4-1}_{l=0}u^{m/2-1-2l}\frac{1}{2^{l}l!}\frac{(m/2)!}{(m/2-2l)!}\frac{J_{(m-2l-3)/2}(v)}{v^{(m-2l-3)/2}},\end{eqnarray*}
 \begin{eqnarray*}C_{m}(u,v)=-\sum^{m/4-1}_{l=0}u^{m/2-1-2l}\frac{1}{2^{l}l!}\frac{(m/2)!}{(m/2-2l)!}\frac{J_{(m-2l-1)/2}(v)}{v^{(m-2l-1)/2}}.\end{eqnarray*}
\end{theorem}
Similarly, we can get the result when $m/2$ is odd.

We can also obtain an alternative expression using exponentials.  When $m$ is even, we have found that formula (\ref{f5}) became   \[\mathcal{L}(t^{m/2-1}e^{i\frac{\pi}{2}\Gamma_{y}}e^{-ti(x,y)})=\frac{\mbox{polynomial of $s$}}{\mbox{polynomial of $s$}}.\]
Hence we can use partial fractions to transform back, as \[\mathcal{L}(t^{m/2-1}e^{i\frac{\pi}{2}\Gamma_{y}}e^{-t(ix,y)})=\sum_{j=1}^{2}\sum_{k=1}^{m/2}\frac{C_{jk}}{(s-\alpha_{j})^{k}}+yx\sum_{j=1}^{2}\sum_{k=1}^{m/2}\frac{C_{jk}}{(s-\alpha_{j})^{k}}.\]
Each $C_{jk}$ can be obtained by the usual technique of partial fractions.

In particular, the kernel of the $2$-dimensional Clifford-Fourier transform can be obtained as follows. Formula (\ref{f5}) becomes \begin{eqnarray*}&&\frac{1}{2\sqrt{+}}\frac{2(s-yx)\sqrt{+}-2\sqrt{+}(x,y)}{s^{2}-(x\wedge y)^{2}}\\&=&\frac{s-yx-(x,y)}{s^{2}-(x\wedge y)^{2}}=\frac{s+(x\wedge y)}{s^{2}-(x\wedge y)^{2}}=\frac{1}{s-(x\wedge y)}.\end{eqnarray*} Transforming back, using (\ref{l7}), we  get the kernel \[K_{2}(x,y)=e^{x\wedge y}.\] This should be compared with section 4.2 in \cite{FNF} and Proposition 5.1 in \cite{HY}.

\subsection{A new integral representation for the Clifford-Fourier kernel}
 When $p=\pi/2$, Theorem \ref{qq1} becomes
\begin{eqnarray}\label{-ix}&&\mathcal{L}(t^{m/2-1}e^{i\frac{\pi}{2}\Gamma_{y}}e^{-ti(x,y)})\nonumber\\&=&\frac{(m/2-1)!}{2\sqrt{+}}\biggl(\frac{s+\sqrt{+}-yx}{(\sqrt{+}+(x,y))^{m/2}}-e^{im\pi/2}\frac{s-\sqrt{+}-yx}{(\sqrt{+}-(x,y))^{m/2}}\biggr)
\end{eqnarray}
By (\ref{l4}), (\ref{l5}) and (\ref{l6}), we have
\begin{eqnarray*} &&\mathcal{L}^{-1}(\frac{s}{\sqrt{+}}\frac{1}{(\sqrt{+}+(x,y))^{m/2}})=\frac{t^{m/2-1}}{(m/2-1)!}e^{-(x,y)t}\\&&-|x||y|t\int_{0}^{t}(t^{2}-u^{2})^{-1/2}J_{1}[|x||y|(t^{2}-u^{2})^{1/2}]\frac{u^{m/2-1}}{(m/2-1)!}e^{-(x,y)u}du
,\end{eqnarray*}
as well as
\begin{eqnarray*} &&\mathcal{L}^{-1}(\frac{1}{\sqrt{+}}\frac{1}{(\sqrt{+}+(x,y))^{m/2-1}})\\&=&\int_{0}^{t}J_{0}[|x||y|(t^{2}-u^{2})^{1/2}]\frac{u^{m/2-2}}{(m/2-2)!}e^{-(x,y)u}du
,\end{eqnarray*}
and
\begin{eqnarray*}&&\mathcal{L}^{-1}(\frac{1}{\sqrt{+}}\frac{x\wedge y}{(\sqrt{+}+(x,y))^{m/2}})\\&=&\int_{0}^{t}(x\wedge y)J_{0}[|x||y|(t^{2}-u^{2})^{1/2}]\frac{u^{m/2-1}}{(m/2-1)!}e^{-(x,y)u}du
.\end{eqnarray*}
Using the above three formulas, we can find the result in the time domain. Then setting $t=1$, we get a new representation of Clifford-Fourier kernel for both even and odd dimension.
\begin{theorem}
The kernel for the $m$-dimensional Clifford-Fourier transform  $(m\ge 3)$ is given by
\begin{eqnarray*}
K_{m}(x,y)=&&\frac{e^{-(x,y)}}{2}-\frac{|x||y|}{2}\int_{0}^{1}(1-u^{2})^{-1/2}J_{1}[|x||y|(1-u^{2})^{1/2}]u^{m/2-1}e^{-(x,y)u}du\\&&+\frac{m-2}{4}\int_{0}^{1}J_{0}[|x||y|(1-u^{2})^{1/2}]u^{m/2-2}e^{-(x,y)u}du
\\&&+\frac{1}{2}\int_{0}^{1}(x\wedge y)J_{0}[|x||y|(1-u^{2})^{1/2}]u^{m/2-1}e^{-(x,y)u}du\\&&-e^{im\pi /2}\biggl(\frac{e^{(x,y)}}{2}-\frac{|x||y|}{2}\int_{0}^{1}(1-u^{2})^{-1/2}J_{1}[|x||y|(1-u^{2})^{1/2}]u^{m/2-1}e^{(x,y)u}du\\&&-\frac{m-2}{4}\int_{0}^{1}J_{0}[|x||y|(1-u^{2})^{1/2}]u^{m/2-2}e^{(x,y)u}du
\\&&+\frac{1}{2}\int_{0}^{1}(x\wedge y)J_{0}[|x||y|(1-u^{2})^{1/2}]u^{m/2-1}e^{(x,y)u}du \biggr).
\end{eqnarray*}
\end{theorem}
\begin{remark} The 2-dimensional kernel was given in the previous section. In this integral representation, the integral is divergent when $m=2$.
\end{remark}
\subsection{Generating function for the kernels of even dimensions}
In this section we  compute the formal generating function of all even dimensional kernels \[G_{p}(x,y,a)=\sum_{m=2,4,6,\cdots}\frac{K^{p}_{m}(x,y)a^{m/2-1}}{\Gamma(m/2)},\]
where $K_{m}^{p}(x,y)$ is the kernel of dimension $m$. Here the formal generating function means one can obtain the kernel from  the derivatives of the generating function. Note that the kernel $K^{p}_{m}(x,y)$ is  in fact a function of $(x,y),$ $|x||y|$ and $|x \wedge y|$. Recall that $\frac{x\wedge y}{|x\wedge y|}$ can be considered as an imaginary unit. So the sum  $G_{p}(x,y,a)$ is not a sum of functions from different spaces but a sum of functions defined on $R^{3}$.

When $p=\pi/2$, the Laplace transform can be computed by
\begin{eqnarray} \label{gn}
&&\sum_{m=2,4,6,\cdots}\frac{1}{(m/2-1)!}\mathcal{L}(t^{m/2-1}e^{i\frac{\pi}{2}\Gamma_{y}}e^{-it(x,y)}a^{m/2-1}) \nonumber\\&=&\frac{1}{2\sqrt{+}}\sum_{m=2,4,6,\cdots}a^{m/2-1}\biggl(\frac{s+\sqrt{+}-yx}{(\sqrt{+}+(x,y))^{m/2}}-e^{im\pi/2}\frac{s-\sqrt{+}-yx}{(\sqrt{+}-(x,y))^{m/2}}\biggr)\nonumber\\&=&\frac{s+\sqrt{+}-yx}{2\sqrt{+}(\sqrt{+}+(x,y)-a)}-e^{i\pi}\frac{s-\sqrt{+}-yx}{2\sqrt{+}(\sqrt{+}-(x,y)-ae^{i\pi})}\nonumber\\&=& \frac{s+\sqrt{+}-yx}{2\sqrt{+}(\sqrt{+}+(x,y)-a)}+\frac{s-\sqrt{+}-yx}{2\sqrt{+}(\sqrt{+}-(x,y)+a)},
\end{eqnarray}
the first equality is by (\ref{-ix}).
Using (\ref{l4}), (\ref{l5}) and (\ref{l6}),  we then get
\begin{eqnarray}\label{f6}&&\mathcal{L}^{-1}\biggl(\frac{s+\sqrt{+}-yx}{\sqrt{+}(\sqrt{+}+(x,y)-a)}\biggr)\nonumber\\&=&2e^{(-(x,y)+a)t}-|x||y|t\int_{0}^{t}(t^{2}-u^{2})^{-1/2} J_{1}[|x||y|(t^{2}-u^{2})^{1/2}]e^{(-(x,y)+a)u}du\nonumber\\&-&|x||y|\int_{0}^{t}e^{(-(x,y)+a)(t^{2}-u^{2})^{1/2}}J_{1}(|x||y|u)du\nonumber\\&-&yx\int_{0}^{t}J_{0}[|x||y|(t^{2}-u^{2})^{1/2}]e^{(-(x,y)+a)u}du
\end{eqnarray}
as well as
\begin{eqnarray}\label{f7}&&\mathcal{L}^{-1}\biggl(\frac{s-\sqrt{+}-yx}{\sqrt{+}(\sqrt{+}-(x,y)+a)}\biggr)\nonumber\\&=&-|x||y|t\int_{0}^{t}(t^{2}-u^{2})^{-1/2} J_{1}[|x||y|(t^{2}-u^{2})^{1/2}]e^{((x,y)-a)u}du\nonumber\\&&+|x||y|\int_{0}^{t}e^{((x,y)-a)(t^{2}-u^{2})^{1/2}}J_{1}(|x||y|u)du\nonumber\\&&+yx\int_{0}^{t}J_{0}[|x||y|(t^{2}-u^{2})^{1/2}]e^{((x,y)-a)u}du
.\end{eqnarray}
Combining (\ref{f6}) and (\ref{f7}), the generating function is
\begin{eqnarray*}G_{\pi/2}(x,y,a)&=&e^{(-(x,y)+a)}\nonumber\\&&-|x||y|\int_{0}^{1}(1-u^{2})^{-1/2} J_{1}[|x||y|(1^{2}-u^{2})^{1/2}] \cosh(((x,y)-a)u)du\nonumber\\&&+|x||y|\int_{0}^{1}\sinh[((x,y)-a)(1-u^{2})^{1/2}]J_{1}(|x||y|u)du\nonumber\\&&+yx\int_{0}^{1}J_{0}[|x||y|(1-u^{2})^{1/2}]\sinh[((x,y)-a)u]du
\end{eqnarray*}
which only gives an integral representation. In the following, we will use  different inverse transform techniques to get the closed form.
Simplifying (\ref{gn}) further,  we have
\begin{eqnarray*}&&\sum_{m=2,4,6,\cdots}\frac{1}{(m/2-1)!}\mathcal{L}(t^{m/2-1}e^{i\frac{\pi}{2}\Gamma_{y}}e^{-it(x,y)}a^{m/2-1})\\&=& \frac{s+\sqrt{+}-yx}{2\sqrt{+}(\sqrt{+}+(x,y)-a)}+\frac{s-\sqrt{+}-yx}{2\sqrt{+}(\sqrt{+}-(x,y)+a)}\\&=&\frac{s-yx-(x,y)+a}{s^{2}+|x|^{2}|y|^{2}-((x,y)-a)^{2}}
.\end{eqnarray*}
Transforming back, we get \begin{eqnarray*}&&\mathcal{L}^{-1}\biggl(\frac{s-yx-(x,y)+a}{s^{2}+|x|^{2}|y|^{2}-((x,y)-a)^{2}}\biggr)\\&=&\cos(\sqrt{|x|^{2}|y|^{2}-((x,y)-a)^{2}}t)\\&&+\frac{-yx-(x,y)+a}{\sqrt{|x|^{2}|y|^{2}-((x,y)-a)^{2}}}\sin(\sqrt{|x|^{2}|y|^{2}-((x,y)-a)^{2}}t)
.\end{eqnarray*}
The last equality is by (\ref{l9}) and  (\ref{l10}).
Note that it equals the case $m=2$ when $a=0$.

Alternatively, a tedious computation shows that \begin{eqnarray*}&&\frac{s-yx-(x,y)+a}{s^{2}+|x|^{2}|y|^{2}-((x,y)-a)^{2}}=\begin{pmatrix}0&1\end{pmatrix}\biggl(s+A\biggr)^{-1}\begin{pmatrix}1\\1\end{pmatrix} .\end{eqnarray*}
Here we introduced the matrix \[A=\biggl((a-(x,y))\begin{pmatrix}1&0\\0&-1\end{pmatrix}+\begin{pmatrix}0 &-xy\\yx& 0\end{pmatrix}\biggr).\]  We get the following
\begin{theorem}\label{th7} The generating function for even dimensional Clifford-Fourier kernels for $p=\frac{\pi}{2}$ is given by
\begin{eqnarray*}G_{\pi/2}(x,y,a)&=&\begin{pmatrix}0&1\end{pmatrix}e^{-A}\begin{pmatrix}1\\1\end{pmatrix}\\&=&\cos(\sqrt{|x|^{2}|y|^{2}-((x,y)-a)^{2}})\\&&+(-yx-(x,y)+a)\frac{\sin(\sqrt{|x|^{2}|y|^{2}-((x,y)-a)^{2}})}{\sqrt{|x|^{2}|y|^{2}-((x,y)-a)^{2}}}
.\end{eqnarray*}
\end{theorem}

We can get a similar result for the fractional case,  i.e. general $p$. Denote \[G_{p}(x,y,a)=\sum_{m=2,4,6,\cdots}\frac{K_{m}^{p}(x,y)a^{m/2-1}}{(m/2-1)!}.\]
Now
\begin{eqnarray*}&&\sum_{m=2,4,\cdots}\mathcal{L}\biggl(\frac{t^{m/2-1}e^{ip\Gamma_{y}}e^{-it(x,y)}a^{m/2-1}}{(m/2-1)!}\biggr)\\&=&\frac{e^{ip}}{2\sqrt{+}}\biggl(\frac{s+\sqrt{+}-ie^{-ip}yx}{s\cos p +i\sqrt{+}\sin p+i(x,y)-ae^{ip}}\\&&-\frac{s-\sqrt{+}-ie^{-ip}yx}{s\cos p-i\sqrt{+}\sin p+i(x,y)-ae^{ip}}\biggr)\\&=&e^{ip}\frac{(-is-e^{-ip}yx)\sin p+(s\cos p +i(x,y)-a e^{ip})}{(s\cos p+i(x,y)-ae^{ip})^{2}+(\sqrt{+})^{2}\sin^{2}p}\\&=&\frac{s-yx\sin p +i(x,y)e^{ip}-ae^{2ip}}{(s\cos p+i(x,y)-ae^{ip})^{2}+(\sqrt{+})^{2}\sin^{2}p}
,\end{eqnarray*}
transforming back by (\ref{l9}) and (\ref{l10}), we have \begin{eqnarray*}&&\mathcal{L}^{-1}\biggl(\frac{s-yx\sin p +i(x,y)e^{ip}-ae^{2ip}}{(s\cos p+i(x,y)-ae^{ip})^{2}+(\sqrt{+})^{2}\sin^{2}p}\biggr)\\&=& e^{-ct}(\cos(dt)+\frac{(x\wedge y-iae^{ip})\sin p}{d}\sin(dt)) \end{eqnarray*}
with $c=(i(x,y)-ae^{ip})\cos p$, $d=\sin p\sqrt{|x|^{2}|y|^{2}+(i(x,y)-ae^{ip})^{2}}$.

Alternatively,  we have
\begin{eqnarray*} &&\frac{s-yx\sin p +i(x,y)e^{ip}-ae^{2ip}}{(s\cos p+i(x,y)-ae^{ip})^{2}+(\sqrt{+})^{2}\sin^{2}p}\\&=&\begin{pmatrix}0&1\end{pmatrix}\biggl(s+B\biggr)^{-1}\begin{pmatrix}1\\1\end{pmatrix}
,\end{eqnarray*}
where \[B=\begin{pmatrix}-\beta_{+}&0\\-(-yx\sin p +i(x,y)e^{ip}-ae^{2ip}+\beta_{+})&-\beta_{-}\end{pmatrix},\] with $\beta_{\pm}$ the roots of $(s\cos p+i(x,y)-ae^{ip})^{2}+(\sqrt{+})^{2}\sin^{2}p,$ i.e.
\[\beta_{\pm}= (-i(x,y)+ae^{ip})\cos p\pm \sin p (\sqrt{-|x|^{2}|y|^{2}-(i(x,y)-ae^{ip})^{2}}).\]
Again, we can have
\begin{theorem}\label{th8} The generating function for the even dimensional fractional Clifford-Fourier kernels is given by
\begin{eqnarray*}G_{p}(x,y,a)&=&\begin{pmatrix}0&1\end{pmatrix}e^{-B}\begin{pmatrix}1\\1\end{pmatrix}\\&=&e^{-c}(\cos d+\frac{(x\wedge y-iae^{ip})\sin p}{d}\sin d) \end{eqnarray*}
with $c=(i(x,y)-ae^{ip})\cos p$ and $d=\sin p\sqrt{|x|^{2}|y|^{2}+(i(x,y)-ae^{ip})^{2}}$.

\end{theorem}

At the end of this section, we give the kernel for general $p$ when $m=2$. It is also the case when $a=0$ in the generating function. The kernel for dimension $2$ is  hence given by \[K_{2}^{p}(x,y)=e^{-i(x,y)\cos p}e^{x\wedge y \sin p},\]
which coincides with the work in \cite{MH} and \cite{HN}.

\

 \footnotesize{\bf Acknowledgments} \qquad Pan Lian is supported by the scholarship from Chinese Scholarship Council (CSC) under the CSC No. 201406120169.

\end{document}